\documentclass[11pt]{amsart}

\usepackage{mystyle}

\date{\today}  

\author{Emma Hogan}
\address{School of Mathematics and Statistics, University of Canterbury, Christchurch, New Zealand}
\curraddr{Mathematical Institute, University of Oxford, Oxford, United Kingdom}
\email{emma.hogan@maths.ox.ac.uk}

\author{Charles Semple}
\address{School of Mathematics and Statistics, University of Canterbury, Christchurch, New Zealand}
\email{charles.semple@canterbury.ac.nz}

\title[Bicircular and lattice path matroids]{The excluded minors for the intersection of bicircular and lattice path matroids}
\thanks{The second author was supported by the New Zealand Marsden Fund.} 

\begin{document}         

\begin{abstract}
  The classes of bicircular matroids and lattice path matroids are closed under minors. The complete list of excluded minors for the class of lattice path matroids is known, and it has been recently shown that the analogous list for the class of bicircular matroids is finite. In this paper, we establish the complete list of excluded minors for the class of matroids that is the intersection of these two classes. This resolves a recently posed open problem. 
\end{abstract}
\maketitle
 
\section{Introduction}

The class of transversal matroids is a classical class of matroids. However, unlike, for example, the classes of representable and graphic matroids, the class of transversal matroids is not closed under minors. As a result, it is natural to study subclasses of transversal matroids with this property. Two particular such classes are the classes of bicircular matroids and lattice path matroids. The former was introduced in \cite{matthews1977bicircular} and the latter was introduced more recently in \cite{Bonin2003}.
 
The complete list of excluded minors for the class of bicircular matroids remains unknown, however it was recently shown that this list is finite \cite{DeVos2021}. The analogous list for the class of lattice path matroids was established in \cite{Bonin2010}. In this paper, we establish the complete list of excluded minors for the class of matroids that is the intersection of bicircular and lattice path matroids, thereby resolving a problem posed in \cite{sivaraman2022}. Other results related to these results include a characterisation of the $3$-connected bicircular matroids whose duals are also bicircular \cite{sivaraman2022} and a characterisation of the matroids that are bicircular and fundamental transversal \cite{neudauer1998transversal}.

During the writing of this paper, a closely related paper appeared on arXiv \cite{guzman2022lattice}. Also motivated by the same problem posed in \cite{sivaraman2022}, the authors investigate, from the viewpoint of bicircular matroids, the intersection of the classes of bicircular and lattice path matroids. The main results characterise those graphs for which the bicircular matroid of the graph is lattice path and characterise the bicircular matroids that are lattice path by listing the minor-minimal bicircular matroids that are not lattice path. It is worth noting that their list of minor-minimal matroids coincides with our list of the bicircular excluded minors for the class of matroids that are bicircular and lattice path. They end the paper asking for the list of minor-minimal lattice path matroids that are bicircular. The main result of this paper also answers their question. 

The paper is organised as follows. In the remainder of this section, we formally state \autoref{thm:main}, the main result of the paper. The next section contains some necessary preliminaries on transversal, lattice path, and bicircular matroids. The proof of \autoref{thm:main} occupies the rest of the paper. In \autoref{sec:lattice_path_EMs}, we establish the list of lattice path excluded minors for the class of bicircular and lattice path matroids. In \autoref{sec:non_lattice_path_EMs}, we establish the list of excluded minors for the class of bicircular and lattice path matroids that are either bicircular, or neither bicircular nor lattice path, thereby completing the proof of \autoref{thm:main}. Throughout the paper, notation and terminology follows Oxley \cite{oxley2011matroid}, with the exception of the definition of bicircular matroids, which, in line with, for example, \cite{sivaraman2022} and \cite{zaslavsky1982bicircular}, we have modified slightly to define a minor-closed class.

To state our main result, we assume the reader is familiar with bicircular and lattice path matroids. Formal definitions are given in the next section. Let $M + e$ denote the free extension of a matroid $M$ by element $e$. The \textit{truncation} $T(M)$ of a matroid $M$ is obtained by taking the free extension $M + e$ of $M$, and then contracting $e$. The \textit{truncation $T_n(M)$ to rank $n$} of a rank-($n+k$) matroid $M$ is obtained by taking a sequence of $k$ truncations of $M$. 

\begin{figure}[H]
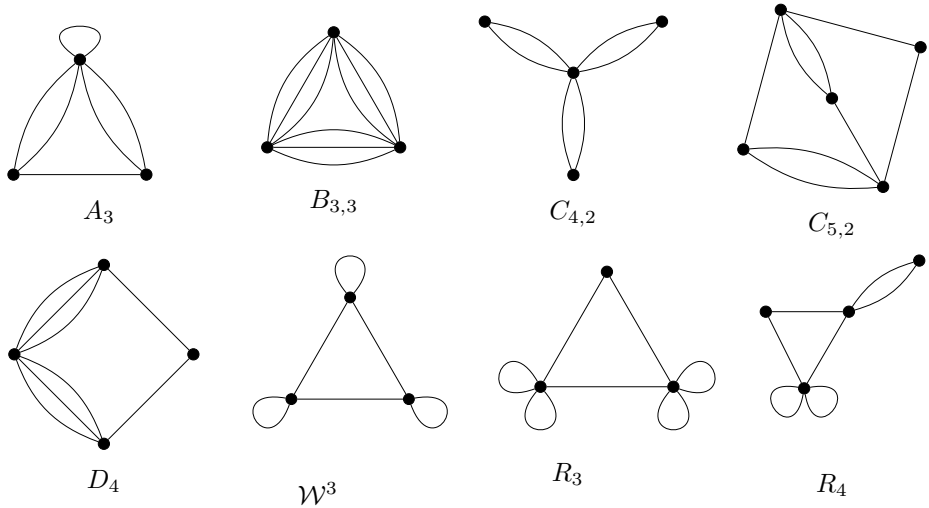

  \centering
  \begin{subfigure}[c]{0.21\textwidth}
    \resizebox{0.85\width}{!}{
      \input{Tikz/ExcludedMinors/a3_bicircular}
    }
    \caption*{$A_3$}
  \end{subfigure}
  \hfill
  \centering
  \begin{subfigure}[c]{0.21\textwidth}
    \centering
    \resizebox{0.85\width}{!}{
      \input{Tikz/ExcludedMinors/b33_bicircular}
    }
    \caption*{$B_{3,3}$}
  \end{subfigure}
  \hfill
  \centering
  \begin{subfigure}[c]{0.22\textwidth}
    \centering
    \resizebox{0.85\width}{!}{
      \input{Tikz/ExcludedMinors/c42_bicircular}
    }
  \caption*{$C_{4,2}$}
  \end{subfigure}
  \hfill
  \centering
  \begin{subfigure}[c]{0.25\textwidth}
    \centering
    \resizebox{0.85\width}{!}{
      \input{Tikz/ExcludedMinors/c52_bicircular}
    }
    \caption*{$C_{5,2}$}
  \end{subfigure}
  \hfill
  \centering
  \begin{subfigure}[c]{0.22\textwidth}
    \centering
    \resizebox{0.85\width}{!}{
      \input{Tikz/ExcludedMinors/d4_bicircular}
    }
    \caption*{$D_4$}
  \end{subfigure}
  \hfill
  \hspace*{-0.6cm}
  \centering
  \begin{subfigure}[c]{0.22\textwidth}
    \centering
    \resizebox{0.82\width}{!}{
      \input{Tikz/ExcludedMinors/whirl3_bicircular}
    }
    \caption*{$\mathcal{W}^3$}
  \end{subfigure}
  \centering
  \hfill
  \begin{subfigure}[c]{0.22\textwidth}
    \centering
    \resizebox{0.85\width}{!}{
      \input{Tikz/ExcludedMinors/r3_bicircular}
    }
    \caption*{$R_3$}
  \end{subfigure}
  \centering
  \hfill
  \begin{subfigure}[c]{0.25\textwidth}
    \centering
    \resizebox{0.85\width}{!}{
      \input{Tikz/ExcludedMinors/r4_bicircular}
    }
    \caption*{$R_4$}
  \end{subfigure}
  \caption{Excluded minors for the class of bicircular and lattice path matroids that are bicircular.}
  \label{fig:excluded_minors_main_2}
\end{figure}

\begin{figure}[H]
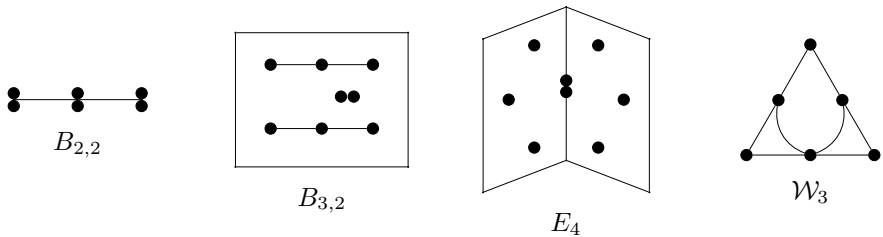

  \centering
  \begin{subfigure}[c]{0.23\textwidth}
      \centering
      \resizebox{0.85\width}{!}{
        \input{Tikz/ExcludedMinors/b22_geometric}
      }
      \caption*{$B_{2,2}$}
  \end{subfigure}
  \hfill
  \centering
  \begin{subfigure}[c]{0.23\textwidth}
      \centering
      \resizebox{0.85\width}{!}{
        \input{Tikz/ExcludedMinors/b32_geometric}
      }
      \caption*{$B_{3,2}$}
  \end{subfigure}
  \hfill
  \centering
  \begin{subfigure}[c]{0.23\textwidth}
      \centering 
      \resizebox{0.85\width}{!}{
        \input{Tikz/ExcludedMinors/e4_geometric}
      }
      \caption*{$E_4$}
  \end{subfigure}
  \hfill
  \centering
  \begin{subfigure}[c]{0.23\textwidth}
      \centering
      \resizebox{0.85\width}{!}{
        \input{Tikz/ExcludedMinors/wheel3_geometric}
      }
    \caption*{$\mathcal{W}_3$}
  \end{subfigure}
  \caption{Excluded minors for the class of bicircular and lattice path matroids that are neither bicircular nor lattice path. Note that $r(B_{3,2})=3$ and $r(E_4) = 4$.}
  \label{fig:excluded_minors_main_3}
\end{figure}

In the statement of \autoref{thm:main}, an excluded minor characterisation of the intersection of bicircular and lattice path matroids, $\mathcal{W}^3$ denotes the rank-$3$ whirl, $\mathcal{W}_3$ denotes the rank-$3$ wheel and $U_{r,n}$ denotes the rank-$r$ uniform matroid on $n$ elements. The remaining named matroids that are not truncations of direct sums of uniform matroids are shown either in \autoref{fig:excluded_minors_main_2} as the bicircular matroid of the given graph, or in \autoref{fig:excluded_minors_main_3} as a geometric representation.

\begin{theorem}
  A matroid is bicircular and lattice path if and only if it has no minor that is isomorphic to any of the matroids
  \begin{enumerate}
    \item[(i)] $U_{3,7}$, $U_{4,7}$, $U_{5,7}$, $T_3(U_{1,2} \oplus U_{3,5})$, $T_3(U_{1,2} \oplus U_{1,2} \oplus U_{3,3})$, $T_4(U_{1,2} \oplus U_{4,5})$, $T_4(U_{3,4} \oplus U_{3,3})$,\\
    \item[(ii)] $A_3$, $B_{3,3}$, $C_{4,2}$, $C_{5,2}$, $D_4$, $\mathcal{W}^3$, $R_3$, $R_4$,\\
    \item[(iii)] $B_{2,2}$, $B_{3,2}$, $E_4$ and $\mathcal{W}_3$.
  \end{enumerate}
  \label{thm:main}
\end{theorem}

We end this section with three remarks. Firstly, the specific matroids listed in \autoref{thm:main} have been grouped so that the matroids in (i) are lattice path, the matroids in (ii) are bicircular, and the matroids in (iii) are neither bicircular nor lattice path. Secondly, to prove \autoref{thm:main}, most of the work is in establishing the lattice path excluded minors. Lastly, multi-path matroids are a generalisation of lattice path matroids first introduced in \cite{Bonin2007multi}, and the class of matroids that are bicircular and multi-path is also suggested in \cite{sivaraman2022} as a potential topic of investigation. Establishing the complete list of excluded minors for this class remains an interesting open problem.

\section{Preliminaries}

The bicircular and lattice path matroids both form minor-closed subclasses of the class of transversal matroids. In this section we define transversal, lattice path, and bicircular matroids, and provide some necessary background on lattice path and bicircular matroids.

Let $E$ be a finite set, and let $J = \{1,2,\dots,r\}$ for some integer $r$. Let $\mathcal{N}$ be a family $(N_j : j \in J)$ of subsets of $E$.
A subset $X = \{x_1,x_2,\dots,x_r\}$ of $E$ is a \textit{transversal} of $\mathcal{N}$ if $x_i \in N_i$ for all $i \in \{1,2,\dots,r\}$. Furthermore, if $K \subseteq J$ and $X$ is a transversal of $(N_j : j \in K)$ then $X$ is a \textit{partial transversal} of $\mathcal{N}$. 
The \textit{transversal matroid} $M[\mathcal{N}]$ is the matroid with ground set $E$ and whose independent sets are the partial transversals of $\mathcal{N}$. We call $\mathcal{N}$ a \textit{presentation} of $M[\mathcal{N}]$.

\noindent \textbf{Lattice path matroids.} Consider a grid, and define a \textit{lattice path} to be a sequence of unit North and East steps on this grid. Let $P$ and $Q$ be lattice paths starting at coordinates $(0,0)$ and 
ending at coordinates $(m, r)$ such that no point on $P$ is strictly above a point on $Q$. Then $P$ and $Q$ define a rank-$r$ matroid on $m+r$ elements as follows.
Consider all the lattice paths $R$ from $(0,0)$ to $(m,r)$ that stay within the region bounded by $P$ and $Q$. For each such lattice path $R$, bijectively label the steps sequentially with the integers from $1$ to $m+r$, and consider the subset of $\{1,2,\ldots,m+r\}$ labelling the North steps of $R$. It is shown in \cite{Bonin2003} that the collection of all such subsets form the bases of a rank-$r$ matroid, $M[P,Q]$. A matroid isomorphic to $M[P,Q]$ is called a lattice path matroid. Furthermore, $(P,Q)$ is the \textit{lattice path presentation} of $M[P,Q]$.

\begin{figure}[H]
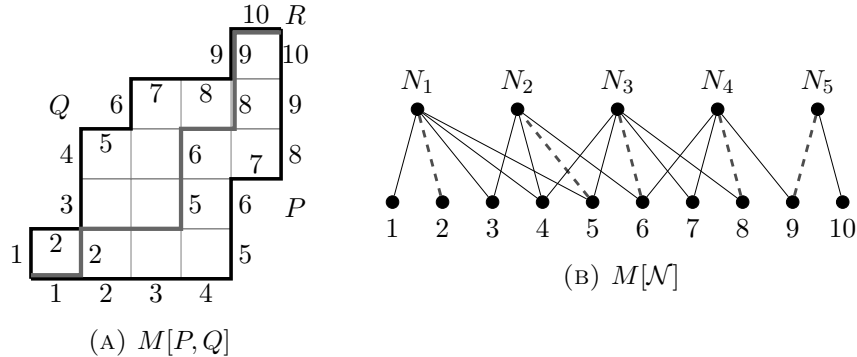

  \centering
    \begin{subfigure}[c]{0.4\textwidth}
      \centering
      \resizebox{.95\width}{!}{
        \input{Tikz/OtherFigures/lattice_path_example.tex}
      }
      \caption{$M[P,Q]$}
      \label{fig:lattice_path_example}
    \end{subfigure}
    \begin{subfigure}[c]{0.55\textwidth}
      \centering
      \resizebox{.95\width}{!}{
        \input{Tikz/OtherFigures/standard_presentation.tex}
      }
      \caption{$M[\mathcal{N}]$}
      \label{fig:standard_presentation}
    \end{subfigure}
    \caption{An example lattice path presentation and the bipartite graph of its standard presentation. 
    Path $R$ represents the basis $\{2,5,6,8,9\}$.}
    \label{fig:lattice_path_presentations}
\end{figure}

It is shown in \cite{Bonin2003} that lattice path matroids are precisely the transversal matroids $M[\mathcal{N}]$ whose ground set admits a linear ordering $x_1,x_2,\dots,x_n$ satisfying two properties. First, each $N_j \in \mathcal{N}$ is the inclusive interval $N_j = [l_j, u_j]$ between some elements $l_j$ and $u_j$ in the ordering of $E$, and second, that the lower and upper endpoints of the intervals $N_j$ form, respectively, the chains $l_1 < l_2 < \dots < l_r$ and $u_1 < u_2 < \dots < u_r$ in the linear ordering. We call $\mathcal{N}$ the \textit{standard presentation} of $M$. This presentation can be obtained from the lattice path presentation $(P,Q)$ by letting $\mathcal{N} = (N_1,N_2,\dots,N_r)$, and taking each $N_i$ to consist of the North steps in the $i$-th row of the region bounded by $P$ and $Q$. In other words, $$N_i = \{x : x \text{ is the $i$-th North step of some path } R\}.$$
An example construction of a standard presentation from the lattice path presentation of a lattice path matroid is given in \autoref{fig:lattice_path_presentations}.

In general, the matroid $M[\mathcal{N}]$ has more than one presentation. However, it is shown in \cite{Bonin2006} that when $M[P,Q]$ is connected, its lattice path presentation $(P,Q)$ is unique up to a $180^\circ$ rotation. In other words, if $P = (p_1,p_2,\dots,p_{m+r})$ and $Q = (q_1,q_2,\dots,q_{m+r})$, then the only other lattice path presentation of $M[P,Q]$ is $(Q^\sigma,P^\sigma)$, where $P^\sigma = (p_{m+r},p_{m+r-1},\dots,p_1)$ and $Q^\sigma = (q_{m+r},q_{m+r-1},\dots,q_1)$. Equivalently, the standard presentation of a connected $M[\mathcal{N}]$ is unique up to reversing the ordering on the ground set, the family $\mathcal{N}$ and each $N_i \in \mathcal{N}$. 

By considering reflections of the lattice path presentation $(P,Q)$ around the $y=x$ axis, it is straightforward to see that the class of lattice path matroids is closed under duality. The following proof that the class of lattice path matroids is also minor-closed is given in \cite{Bonin2006}. We include the proof here as it establishes a certain presentation of a single-element deletion of a lattice path matroid $M$ from a presentation of $M$.

\begin{lemma}
    The class of lattice path matroids is closed under minors. \label{lem:lattice_path_closed_minors}
  \end{lemma}
  \begin{proof}
    Let $M$ be a lattice path matroid with ground set $\{1,2,\dots,m+r\}$ and standard presentation $\cN = (N_1,N_2,\dots,N_r)$ where $N_i = [l_i, u_i]$ for each $i \in \{1,2,\ldots,r\}$. Suppose $x$ is in the ground set of $M$. Since the class of lattice path matroids is closed under duality, it suffices to show that $M\backslash x$ is lattice path.
  
    A presentation of $M \backslash x$ is given by $$\cN' = (N_1',N_2',\dots,N_r') = (N_1 - x, N_2 - x,\dots,N_r - x).$$ This is the standard presentation of a lattice path matroid if and only if both the lower and upper endpoints of $\cN'$ form chains in the induced linear order $[m+r] - x$. In all cases where this is not true, we show that there exists a modification $\cN''$ of $\cN'$ matching the same sets, and for which the lower and upper endpoints do form chains in $[m+r] - x$.
  
    First, if $N_i' = \emptyset$ for some $i \in \{1,2,\ldots,r\}$, then $x$ is a co-loop of $M$ and so $N_i$ is the only interval in $\cN$ containing $x$. Thus, $\cN'' = \cN' - N_i'$ is a standard presentation of $M\backslash x$, and so $M\backslash x$ is lattice path. Otherwise, every interval $N_i'$ in $\cN'$ is non-empty. If the upper endpoints, $u_i'$ of $\cN'$ do not form a chain in $[m+r] - x$, then we must have that $x$ and $x-1$ are respectively the upper endpoints $u_k$ and $u_{k-1}$ for 
    some $k$ in $\{2,3,\ldots,r\}$. 
  
    Let $B$ be a basis of $M$ not containing $x$. Then, since $M$ is lattice path, we may write $B = \{b_1, b_2, \ldots, b_r\}$ so that $b_1 < b_2 < \cdots < b_r$ and $b_i \in N_i$ for each $i$ in $\{1,2,\dots,r\}$. Since $b_k \in N_k$ and $x = u_k$ is not in $B$, we have that $b_k \leq x-1$. Hence, if $x-1$ is in $B$, then $x - 1 = b_q$ for some $q \geq k$, and so in all cases we have that $B$ can be matched with the presentation $$(N_1',N_2',\dots,N_{k-2}', N_{k-1}' - (x-1), N_k',\dots,N_r').$$
    Similarly, if $x-2 = u_{k-2}$, we can replace $N_{k-2}'$ by $N_{k-2}' - (x-2)$. In general, if $k-l$ is the least value for which $u_{k-l} = x - l$, we can take 
    \begin{align*}
        \cN'' &= (N_1',N_2',\dots,N_{k-l-1}',N_{k-l}' - (x - l), \\
        &\qquad \qquad  N_{k-l+1}' - (x - l + 1),\dots, N_{k-1}' - (x-1), N_k', N_{k+1}', \dots, N_r').
    \end{align*} 
    to be a presentation of $M\backslash x$ for which the upper endpoints form a chain in $[m+r] - x$. Symmetrically, modifications to the lower endpoints of intervals in $\cN'$ can be made to obtain a presentation of $M\backslash x$ in which the lower endpoints also form a chain in $[m+r] - x$. Thus, $M\backslash x$ has a standard presentation as a lattice path matroid, and so the class of lattice path matroids is closed under minors.
  \end{proof}

  \noindent \textbf{Remark}. The proof of \autoref{lem:lattice_path_closed_minors} showed how a certain presentation of a single-element deletion of a lattice path matroid $M$ can be obtained from a presentation of $M$. In particular, we have the following:
  
  Let $M$ be a lattice path matroid with standard presentation $\mathcal{N} = (N_1,N_2,\dots,N_r)$, and let $x$ be the upper endpoint of the interval $N_k$. Let $N_i' = N_i - x$ for all $i \in \{1,2,\ldots, r\}$. If $k-l$ is the least value for which $u_{k-l} = x - l$, then 
  \begin{align*}
        \cN'' &= (N_1',N_2',\dots,N_{k-l-1}',N_{k-l}' - (x - l), \\
        &\qquad \qquad  N_{k-l+1}' - (x - l + 1),\dots, N_{k-1}' - (x-1), N_k', N_{k+1}', \dots, N_r')
  \end{align*}
  is a presentation of $M\backslash x$ in which the upper endpoints of $\mathcal{N}'$ form a chain. We make use of this construction in the proof of \autoref{lem:uniform_minor}. \qed

\begin{figure}[H]
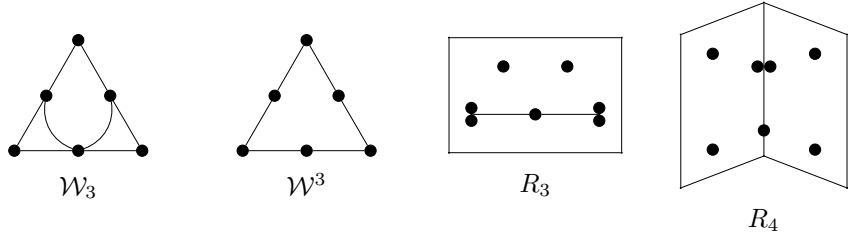

  \centering
  \begin{subfigure}[c]{0.23\textwidth}
    \centering
    \resizebox{.85\width}{!}{
      \input{Tikz/ExcludedMinors/wheel3_geometric}
    }
  \caption*{$\mathcal{W}_3$}
  \end{subfigure}
  \centering
  \begin{subfigure}[c]{0.23\textwidth}
      \centering
      \resizebox{.85\width}{!}{
        \input{Tikz/ExcludedMinors/whirl3_geometric}
      }
    \caption*{$\mathcal{W}^3$}
  \end{subfigure}
  \centering
  \begin{subfigure}[c]{0.23\textwidth}
      \centering
      \resizebox{.85\width}{!}{
        \input{Tikz/ExcludedMinors/r3_geometric}
      }
    \caption*{$R_3$}
  \end{subfigure}
  \centering
  \begin{subfigure}[c]{0.23\textwidth}
      \centering
      \resizebox{.85\width}{!}{
        \input{Tikz/ExcludedMinors/r4_geometric}
      }
    \caption*{$R_4$}
  \end{subfigure}
  \caption{Excluded minors for the class of lattice path matroids that are not in infinite families.}
  \label{fig:individual_L_excluded_minors}
\end{figure}

An excluded-minor characterisation for the class of lattice path matroids is given in \cite{Bonin2010}. 
The characterisation consists of four individual excluded minors and five infinite families of excluded minors.
The individual excluded minors are the rank-$3$ wheel $\mathcal{W}_3$, the rank-$3$ whirl $\mathcal{W}^3$, $R_3$, and $R_4$. Geometric representations of these matroids are shown in \autoref{fig:individual_L_excluded_minors}. Note that $R_4$ is a rank-$4$ matroid, while the remainder are rank $3$.

\begin{figure}[H]
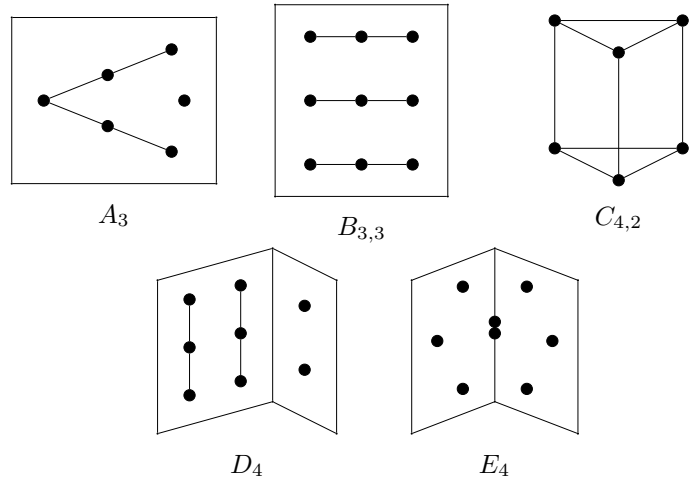

  \centering
  \begin{subfigure}[c]{0.25\textwidth}
    \centering
    \resizebox{.85\width}{!}{
      \input{Tikz/ExcludedMinors/a3_geometric}
    }
  \caption*{$A_3$}
  \end{subfigure}
  \centering
  \begin{subfigure}[c]{0.25\textwidth}
      \centering
      \resizebox{.85\width}{!}{
        \input{Tikz/ExcludedMinors/b33_geometric}
      }
    \caption*{$B_{3,3}$}
  \end{subfigure}
  \centering
  \begin{subfigure}[c]{0.27\textwidth}
      \centering
      \resizebox{.85\width}{!}{
        \input{Tikz/ExcludedMinors/c42_geometric}
      }
    \caption*{$C_{4,2}$}
  \end{subfigure}
  \centering
  \begin{subfigure}[c]{0.25\textwidth}
      \centering
      \resizebox{.85\width}{!}{
        \input{Tikz/ExcludedMinors/d4_geometric}
      }
    \caption*{$D_4$}
  \end{subfigure}
  \begin{subfigure}[c]{0.25\textwidth}
    \centering
    \resizebox{.85\width}{!}{
      \input{Tikz/ExcludedMinors/e4_geometric}
    }
  \caption*{$E_4$}
\end{subfigure}
  \caption{Examples of excluded minors for the class of lattice path matroids from infinite families.}
  \label{fig:families_L_excluded_minors}
\end{figure}

We now define the five infinite families of excluded minors for the class of lattice path matroids, labelled as in \cite{Bonin2010}. Let $P_n$ denote the matroid $T_n(U_{n-1,n} \oplus U_{n-1,n})$, and let $P_n'$ denote the matroid $(P_{n-1}^* + e)^*$, the free single-element co-extension of $P_{n-1}$. The five infinite families of excluded minors for the class of lattice path matroids are as follows:
\begin{enumerate}
  \item $A_n = P_n' + x$ for all $n \geq 3$,
  \item $B_{n,k} = T_n(U_{n-1,n} \oplus U_{n-1,n} \oplus U_{k-1,k})$ for all $n \geq k \geq 2$,
  \item $C_{n+k,k} = B_{n,k}^*$ for all $n \geq k \geq 2$,
  \item $D_n = (P_{n-1} \oplus U_{1,1}) + x$ for all $n \geq 4$,
  \item $E_n = D_n^*$ for all $n \geq 4$.
\end{enumerate}
An example of each of these families, as a geometric representation, is given in \autoref{fig:families_L_excluded_minors}.

The next theorem is the excluded-minor characterisation of the class of lattice path matroids established in \cite{Bonin2010}.
\begin{theorem}
  A matroid $M$ is lattice path if and only if it has no minor isomorphic to any of the matroids $\mathcal{W}_3$, $\mathcal{W}^3$, $R_3$, $R_4$, or any matroid in the families $A_n$ for all $n \geq 3$, $B_{n,k}$ for all $n \geq k \geq 2$, $C_{n+k,k}$ for all $n \geq k \geq 2$, $D_n$ for all $n \geq 4$, and $E_n$ for all $n \geq 4$.
  \label{thm:lattice_path_excluded_minors}\qed
\end{theorem}

\noindent \textbf{Bicircular matroids.} Let $G$ be a graph with edge set $E(G)$. We define the edge set $E(G)$ to include edges of three distinct types: \textit{links}, which are edges with two distinct endpoints, \textit{loops} which are edges with both endpoints at the same vertex, and \textit{free edges} which are edges adjacent to no vertices. Modifying the definition of a graph to include free edges allows us to define the bicircular matroids as a minor-closed class, in contrast to the definition provided in \cite{oxley2011matroid}.

\begin{figure}[H]
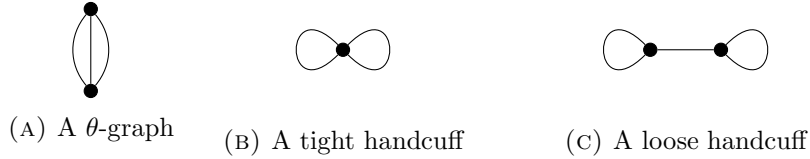

  \centering
  \hspace*{-1cm}
  \begin{subfigure}[c]{.3\textwidth}
      \centering
      \input{Tikz/OtherFigures/theta_graph.tex}
      \caption{A $\theta$-graph}
  \end{subfigure}
  \hspace*{-0.7cm}
  \begin{subfigure}[c]{.3\textwidth}
      \centering
      \input{Tikz/OtherFigures/tight_handcuff.tex}
      \caption{A tight handcuff}
  \end{subfigure}
  \hspace*{0.5cm}
  \begin{subfigure}[c]{.3\textwidth}
      \centering
      \input{Tikz/OtherFigures/loose_handcuff.tex}
      \caption{A loose handcuff}
  \end{subfigure}
  \caption{Bicircular graphs are subdivisions of one of these three minimal examples.}
  \label{fig:bicircular_graphs}
\end{figure}

A \textit{bicircular graph} is a subdivision of any of the graphs in \autoref{fig:bicircular_graphs}. These graphs come in three classes: \textit{$\theta$-graphs}, which are two cycles sharing a common path of at least one edge, \textit{tight handcuffs} which are two cycles sharing a common vertex but no common edge, and \textit{loose handcuffs} which are graphs consisting of two cycles joined by a non-zero length path.

The bicircular matroid of a graph $G$ is the matroid whose circuits are the free edges of $G$ (representing matroid loops), and the subsets of $E(G)$ inducing a bicircular graph. We call $G$ a \textit{bicircular representation} of $B(G)$.

The next two lemmas provide alternate representations of bicircular matroids. In general, neither these representations, nor the bicircular representation of a bicircular matroid are unique. The following result of \cite{matthews1977bicircular} establishes the bicircular matroids as a subclass of the transversal matroids.

\begin{lemma}
  Let $M$ be a matroid. Then $M$ is bicircular if and only if $M$ has a transversal presentation $\mathcal{N}$ such that each element of $M$ is contained in at most two sets of $\mathcal{N}$. \label{lem:bicircular_presentation} 
\end{lemma}

A \textit{simplex} $\Delta$ is a geometric representation of a rank-$r$ matroid $M$ if $\Delta$ consists of $r$ vertices placed in $(r-1)$-dimensional Euclidean space and the elements of each circuit $C$ of $M$ are placed in a $(r(C) - 1)$-dimensional `face' of $\Delta$. The following result is also established in \cite{matthews1977bicircular}.

  \begin{lemma}
      Let $M$ be a rank-$r$ matroid. Then $M$ is bicircular if and only if $M$ can be 
      represented geometrically by placing all non-loop elements of $M$ on the vertices or edges of an $r$-vertex 
      simplex $\Delta$.
      \label{lem:bicircular_simplex}
  \end{lemma} 
  
  The complete set of excluded minors for the class of bicircular matroids is not known. However, \cite{DeVos2021} shows that the size of this set is finite and lists $27$ known excluded minors for the class. We will use the following subset of this list in this paper.

  \begin{lemma}
    The matroids $U_{3,7}$, $U_{5,7}$, $T_3(U_{1,2} \oplus U_{3,5})$, $T_3(U_{1,2} \oplus U_{1,2} \oplus U_{3,3})$, $T_4(U_{1,2} \oplus U_{4,5})$, $T_4(U_{3,4} \oplus U_{3,3})$, $B_{2,2}$, $B_{3,2}$, $E_4$ and $\mathcal{W}_3$ are excluded minors for the class of bicircular matroids.
    \label{lem:bicircular_excluded_minors}
  \end{lemma}

  Let $M$ be a matroid on ground set $E$, and let $k$ be a positive integer. A partition $(A, B)$ of $E(M)$ is a \textit{k-separation} of $M$ if $|A| \geq k$, $|B| \geq k$ and $r(A) + r(B) - r(M) < k$. If additionally, $r(A) \geq k$ and $r(B) \geq k$, we call $(A,B)$ a 
  \textit{vertical k-separation}. We say $M$ is \textit{vertically k-connected} if $M$ has no vertical $l$-separations for $l < k$.
  The following lemma about excluded minors for the class of bicircular matroids is established in \cite{DeVos2021}.

  \begin{lemma} \label{lem:parallel_vert_3}
    Let $M$ be an excluded minor for the class of bicircular matroids. Then 
    \begin{enumerate}
      \item[\rm{(i)}] Every parallel class of $M$ consists of at most two elements, and 
      \item[\rm{(ii)}] $M$ is vertically $3$-connected. 
    \end{enumerate}
  \end{lemma}

\section{Excluded Minors that are Lattice Path} \label{sec:lattice_path_EMs}

In this section we show that the lattice path excluded minors for the class of bicircular matroids are exactly $U_{3,7}$, $U_{4,7}$, $U_{5,7}$, $T_3(U_{1,2}\oplus U_{3,5})$, $T_3(U_{1,2} \oplus U_{1,2} \oplus U_{3,3})$, $T_4(U_{1,2} \oplus U_{4,5})$ and $T_4(U_{3,4} \oplus U_{3,3})$. First note that, by \autoref{lem:bicircular_excluded_minors}, these matroids are excluded minors for the class of bicircular matroids. Furthermore, lattice path presentations of each of these matroids are given in \autoref{fig:lattice_path_excluded_minors}, along with geometric representations of those that are non-uniform. This gives us the following result:
\begin{lemma}
  The matroids $U_{3,7}$, $U_{4,7}$, $U_{5,7}$, $T_3(U_{1,2}\oplus U_{3,5})$, $T_3(U_{1,2} \oplus U_{1,2} \oplus U_{3,3})$, $T_4(U_{1,2} \oplus U_{4,5})$ and $T_4(U_{3,4} \oplus U_{3,3})$ are all lattice path excluded minors for the class of bicircular and lattice path matroids. \label{lem:lattice_path_excluded_minors}
  \qed
\end{lemma}

We spend the remainder of this section proving that the list given in \autoref{lem:lattice_path_excluded_minors} is the complete list of such matroids. It is shown in \cite{Bonin2006} and \cite{matthews1977bicircular}, respectively, that the classes of lattice path and bicircular matroids are both closed under direct sums, and so the intersection of these two classes is also closed under direct sums. Therefore, all excluded minors for the class of bicircular and lattice path matroids are connected. It follows that the lattice path presentations of excluded minors in this section are all unique up to a $180^\circ$ rotation.

\begin{figure}[H]
  \centering 
  \begin{subfigure}[c]{0.3\textwidth}
    \centering
    \resizebox{0.9\width}{!}{
      \input{Tikz/ExcludedMinors/lattice_path1}
    }
    \caption*{$U_{3,7}$}
  \end{subfigure}
  \hfill
  \begin{subfigure}[c]{0.3\textwidth}
    \centering
    \resizebox{0.9\width}{!}{
      \input{Tikz/ExcludedMinors/lattice_path4}
    }
    \caption*{$U_{4,7}$}
  \end{subfigure}
  \hfill
  \begin{subfigure}[c]{0.3\textwidth}
    \centering
    \resizebox{0.9\width}{!}{
      \input{Tikz/ExcludedMinors/lattice_path7}
    }
    \caption*{$U_{5,7}$}
    \vspace{0.2cm}
  \end{subfigure}

  \centering
  \begin{subfigure}[c]{0.45\textwidth}
    \centering
    \begin{subfigure}[c]{0.45\textwidth}
      \centering
      \resizebox{0.9\width}{!}{
        \input{Tikz/ExcludedMinors/lattice_path2}
      }
    \end{subfigure}
    \hfill
    \begin{subfigure}[c]{0.45\textwidth}
      \centering
      \resizebox{0.8\width}{!}{
        \input{Tikz/ExcludedMinors/geometric2}
      }
    \end{subfigure}
    \caption*{$T_3(U_{1,2} \oplus U_{3,5})$}
  \end{subfigure}
  \hfill
  \begin{subfigure}[c]{0.45\textwidth}
    \centering
    \begin{subfigure}[c]{0.45\textwidth}
      \centering
      \resizebox{0.9\width}{!}{
        \input{Tikz/ExcludedMinors/lattice_path3}
      }
    \end{subfigure}
    \hfill
    \begin{subfigure}[c]{0.45\textwidth}
      \centering
      \resizebox{.8\width}{!}{
        \input{Tikz/ExcludedMinors/geometric3}
      }
    \end{subfigure}
    \caption*{$T_3(U_{1,2} \oplus U_{1,2} \oplus U_{3,3})$}
  \end{subfigure}

  \centering
  \begin{subfigure}[c]{0.45\textwidth}
    \vspace{0.3cm}
    \centering
    \begin{subfigure}[c]{0.45\textwidth}
      \centering
      \resizebox{0.9\width}{!}{
        \input{Tikz/ExcludedMinors/lattice_path6}
      }
    \end{subfigure}
    \begin{subfigure}[c]{0.45\textwidth}
      \centering
      \resizebox{.8\width}{!}{
        \input{Tikz/ExcludedMinors/geometric6_v2}
      }
    \end{subfigure}
    \caption*{$T_4(U_{1,2} \oplus U_{4,5})$}
  \end{subfigure}
  \hfill
  \begin{subfigure}[c]{0.45\textwidth}
    \vspace{0.4cm}
    \centering
    \begin{subfigure}[c]{0.45\textwidth}
      \centering
      \resizebox{0.9\width}{!}{
        \input{Tikz/ExcludedMinors/lattice_path5}
      }
    \end{subfigure}
    \begin{subfigure}[c]{0.45\textwidth}
      \centering
      \resizebox{0.8\width}{!}{
        \input{Tikz/ExcludedMinors/geometric5}
      }
    \end{subfigure} 
    \caption*{$T_4(U_{3,4} \oplus U_{3,3})$}
  \end{subfigure}
  \caption{Excluded minors for the class of bicircular and lattice path matroids that are lattice path.}
  \label{fig:lattice_path_excluded_minors}
\end{figure}

We begin by showing that $U_{5,7}$ is the only lattice path excluded minor for the intersection of bicircular and lattice path matroids with rank at least $5$. This involves a sequence of lemmas culminating in \autoref{lem:only_u_57}.

\begin{lemma}
    Let $M$ be a rank-$r$ lattice path matroid with standard presentation $\mathcal{N} = (N_1,N_2,\dots,N_r)$. If $M$ is not bicircular, then there exists some element $x$ in $E(M)$ such that $x \in N_i \cap N_{i+1} \cap N_{i+2}$ for some $i \in \{1,2,\dots,r-2\}$. 
    \label{lem:2x3_subgrid}
  \end{lemma}
  \begin{proof}
    If each element of $M$ is contained in at most two intervals of $\mathcal{N}$, then \autoref{lem:bicircular_presentation} implies that $M$ is bicircular. Therefore, if $M$ is not bicircular, there must exist at least one element $x$ contained in at least three intervals of $\mathcal{N}$. Moreover, the intervals containing $x$ are consecutive in $\mathcal{N}$.
  \end{proof}
  
  Let $M$ be a rank-$r$ lattice path matroid with ground set $\{1, 2, \ldots, r+m\}$ and standard presentation $\cN=(N_1, N_2, \ldots, N_r)$. We say that $\cN$ has the {\em upper bound property} if $l_{k+2} \leq u_k$ for all $k \in \{1, 2, \ldots, r-2\}$.
  
  \begin{lemma}
      Let $M$ be a rank-$r$ lattice path matroid with ground set $\{1, 2, \ldots, r+m\}$ and standard presentation $\cN = (N_1, N_2, \ldots, N_r)$. If $M$ is vertically 3-connected, then $\cN$ has the upper bound property.
      \label{lem:vert_3_connected_interval}
  \end{lemma}
          
  \begin{proof}
      Let $M$ be vertically $3$-connected, and suppose that there exists a $k\in \{1, 2, \ldots, r-2\}$ such that $l_{k+2} > u_k$. Let $A=\{1, 2, \ldots, l_{k+2}-1\}$ and $B=\{l_{k+2}, l_{k+2}+1, \ldots, r+m\}$, 
      so that $(A, B)$ is a partition of $E(M)$. We have $l_{k+1} \leq l_{k+2} - 1$ by definition and so $\{l_1,l_2,\ldots,l_{k+1}\} \subseteq A$. Furthermore, since $l_{k+2} > l_{k+2} - 1$, the set $A$ does not meet $N_{k+2}$ and hence $r(A) = k+1$. If $l_{k+2} > u_{k+1}$ then $M$ is disconnected, a contradiction, so $\{u_{k+1}, u_{k+2} \ldots, u_r\} \subseteq B$. Since $l_{k+2} > u_k$, the set $B$ does not meet $N_{k}$ and hence $r(B) = r-k$. Therefore, $r(A)+r(B)-r(M)=1$.
  
      Finally, as $k\in \{1, 2, \ldots, r-2\}$, we have $r(A) \geq 2$ and $r(B) \geq 2$. Hence, $(A, B)$ is a vertical $2$-separation, a contradiction. This completes the proof 
      of the lemma.
  \end{proof}
      
  \begin{lemma}
    Let $M$ be a rank-$r$ lattice path matroid with ground set $\{1, 2, \ldots, r+m\}$ and standard presentation $\cN = (N_1, N_2, \ldots, N_r)$. If $r > 2$ and $\cN$ has the upper bound property, then $|N_i|\geq 3$ for all $i\in \{1, 2, \ldots, r\}$.
    \label{lem:vert_3_connected_interval_size}
  \end{lemma}
      
  \begin{proof}
      Suppose that, for some $i \in \{1, 2, \ldots, r\}$, we have $|N_i| < 3$. Then $u_i < l_i+2$. If $i \in \{1, 2, \ldots, r-2\}$, then $l_i+2 \leq l_{i+2}$, and so 
      $u_i < l_{i+2}$, contradicting the upper bound property. Thus, we may assume that $i \in \{r-1, r\}$.
      
      If $i=r-1$, then $u_{r-1} < l_{r-1}+2$. Therefore, as $l_{r-1} + 1 \leq l_r$ and $u_{r-2} \leq u_{r-1} - 1$, we have
      $$u_{r-2} \leq u_{r-1} - 1 < l_{r-1}+1 \leq l_r.$$
      Thus, $u_{r-2} < l_r$, contradicting the upper bound property. Lastly, if $i=r$, then $u_r < l_r+2$. But $u_{r-2} \leq u_r-2$, and so
      $$u_{r-2} \leq u_r - 2 < l_r$$
      that is, $u_{r-2} < l_r$. This last contradiction to the upper bound property completes the proof of the lemma.
  \end{proof}
  \newpage
  \begin{lemma}
    Let $M$ be a rank-$r$ lattice path matroid with $r \geq 3$. If $M$ is vertically $3$-connected, then $M$ has a $U_{r, r+2}$-restriction.
    \label{lem:uniform_minor}
  \end{lemma}
        
  \begin{proof}
    Let $\{1, 2, \ldots, r+m\}$ denote the ground set of $M$, and let $\cN=(N_1, N_2, \ldots, N_r)$ be a standard presentation of $M$, where $N_i=[l_i, u_i]$ for all $i\in \{1, 2, \ldots, r\}$. Suppose that $M$ is vertically $3$-connected. Then \autoref{lem:vert_3_connected_interval} implies that $\cN$ has the upper bound property. Consider the lower endpoints of the intervals $N_1, N_2, \ldots, N_r$, and assume that there exists an $i\in \{1, 2, \ldots, r-1\}$ such that $l_{i+1}\neq l_i+1$. Let $k$ be the smallest index for which $l_{k+1}\neq l_k+1$. We next show that the rank-$r$ lattice path matroid $M\backslash (l_k+1)$ has a standard presentation satisfying the upper bound property.
    
    Since $l_k < l_k+1 < l_{k+1}$, and the lower endpoints of $N_1, N_2, \ldots, N_r$ form a chain, $l_k+1$ is not the lower endpoint of any interval in $\cN$. If $l_k+1$ is not the upper endpoint of any interval in $\cN$, then
    \begin{align*}
      \cN' &= (N_1-(l_k+1), N_2-(l_k+1), \ldots, N_r-(l_k+1))
    \end{align*}
    is a standard presentation of $M\backslash (l_k+1)$. Since $\cN$ satisfies the upper bound property, it follows that $(N_1-(l_k+1), N_2-(l_k+1), \ldots, N_r-(l_k+1))$ also satisfies this property under the induced linear order on $[m+r]-(l_k+1)$.
    
    On the other hand, suppose $l_{k} + 1$ is the upper endpoint of an interval in $\cN$. Say $l_{k} + 1 = u_j$. Since $l_k+1 < l_{k+1}$ and $l_{k+1} \leq u_{k-1}$ by the upper bound property, it follows that $u_j < u_{k-1}$, and so $j \leq k-2$. Since $u_j = l_k + 1$ is only an element of the intervals $N_j, N_{j+1},\ldots,N_k$, \autoref{lem:lattice_path_closed_minors} implies that a standard presentation of $M\backslash (l_k+1)$ is given by
    \begin{align*}
    \cN' & = (N'_1, N'_2, \ldots, N'_r) \\
    & = (N_1, N_2, \ldots, N_{j-c-1}, N_{j-c}-u_{j-c}, \ldots, N_{j-1}-u_{j-1}, \\
    & \qquad \qquad \qquad N_j-u_j, N_{j+1}-u_j, \ldots, N_k-u_j, N_{k+1}, N_{k+2}, \ldots, N_r),
    \end{align*}
    where $c \geq 0$ and $j-c$ is the least value such that $u_{j-c}=u_j-c$. Let $l'_i$ and $u'_i$ denote the lower and upper endpoints, respectively, of $N'_i$ for all $i\in \{1, 2, \ldots, r\}$. To establish that $\cN'$ satisfies the upper bound property (under the induced linear order), it is easily checked that it suffices to show that $l'_{i+2}\le u'_i$ for all $i\in \{j-c, j-c+1, \ldots, j\}$.
    
    Since $j\le k-2$, it follows that $l'_{j+2}=l_{j+2}\le l_k=u_j-1=u'_j$, and so $l'_{j+2}\le u'_j$. For all $i\le j-1$, let $i=j-a$, where $a\ge 1$. Then $u_i=u_j-a$ and $u'_i=u_j-a-1=u'_j-a$. Similarly, $l'_{j+2-a}=l'_{j+2}-a$. Thus, $l'_{j+2}\le u'_j$ implies that $l'_{j+2}-a\le u'_j-a$, and so $l'_{i+2}\le u'_i$ for all $i\in \{j-c, j-c+1, \ldots, j\}$.
    
    Applying this process to the next element in $[m+r]-(l_k+1)$ that is not a lower endpoint of an interval in $\cN'$ and repeating until there are no such elements, we eventually obtain a rank-$r$ restriction $M'$ of $M$ with a standard presentation whose lower endpoints are consecutive in the induced linear order on the resulting ground set. Observe that if, for some $i$, we have $u_{i+1}=u_i+1$ in $\cN$, then the upper endpoints of the $i$-th and $(i+1)$-th intervals of the standard presentation of $M'$ remain consecutive in the induced linear order. In other words, we can symmetrically apply the same argument to the upper endpoints of this standard presentation of $M'$ without changing the consecutive ordering of the lower endpoints that we achieved in the standard presentation of $M'$. Hence, we obtain a restriction $M''$ of $M'$, and thus of $M$, with a standard presentation $(N''_1, N''_2, \ldots, N''_r)$ whose upper endpoints, as well as lower endpoints, are consecutive in the induced linear order on the resulting ground set, and which satisfies the upper bound property. Therefore,
    $$|N''_1|=|N''_2|=\cdots =|N''_r|.$$
    By \autoref{lem:vert_3_connected_interval_size}, $|N''_i|\ge 3$ for all $i\in \{1, 2, \ldots, r\}$, and so $M''$ is isomorphic to $U_{r, r+t}$ for some $t\ge 2$. Hence, $M$ has a $U_{r, r+2}$-restriction.
  \end{proof}
  
  \begin{lemma}
    Let $M$ be a rank-$r$ lattice path matroid with $r \geq 3$. If $M$ is an excluded minor for the class of bicircular matroids then $M$ has a $U_{r, r+2}$-restriction.
    \label{lem:uniform_minor_bicircular}
  \end{lemma}
  \begin{proof}
    \autoref{lem:parallel_vert_3} implies that $M$ is vertically $3$-connected. Hence, the result follows from \autoref{lem:uniform_minor}.
  \end{proof}
  
  \begin{lemma}
    Let $M$ be a rank-$r$ lattice path excluded minor for the class of bicircular and lattice path matroids. If $r \geq 5$, then $M$ is isomorphic to $U_{5,7}$.
    \label{lem:only_u_57}
  \end{lemma}
  \proof 
  First, note that $U_{5,7}$ is lattice path, with the lattice path presentation given in \autoref{fig:lattice_path_excluded_minors}, so by \autoref{lem:bicircular_excluded_minors}, it is an excluded minor for the class of bicircular and lattice path matroids. Now suppose $M$ is a rank-$r$ lattice path excluded minor for this class, and let $r \geq 5$. Then \autoref{lem:uniform_minor_bicircular}
  implies that $M$ has a $U_{r,r+2}$-restriction, and thus a $U_{5,7}$-minor. Since $U_{5,7}$ is itself an excluded minor for the class of bicircular and lattice path matroids, we must have $M \cong U_{5,7}$. Hence, $U_{5,7}$ is the only lattice path matroid of rank greater than $4$ that is an excluded minor for the class of bicircular and lattice path matroids. \qed
  
  We now show that we have a complete set of the lattice path excluded minors for ranks less than $5$.

  \begin{lemma}
    Let $M$ be a rank-$r$ lattice path matroid with standard presentation $\mathcal{N} = (N_1,N_2,\dots,N_r)$, where $N_i = [l_i, u_i]$ for all $i \in \{1,2,\dots,r\}$. Suppose $M$ is an excluded minor for the class of bicircular matroids. If $\{a,b\}$ is a parallel class of $M$, then either $\{a,b\} = \{l_1,l_1+1\}$, in which case $l_2 = l_1 + 2$, or $\{a,b\} = \{u_r, u_r - 1\}$, in which case $u_r = u_{r-1} + 2$.
    \label{lem:gross_parallel_class}
  \end{lemma}
  \begin{proof}
    Since $M$ is an excluded minor for the class of bicircular matroids, \autoref{lem:parallel_vert_3} implies that any parallel class of $M$ contains at most two elements. Suppose $\{a,b\}$ is a parallel class of $M$ and assume $a < b$. Then both $a$ and $b$ are members of exactly one interval $N_i$ in $\mathcal{N}$. Because the parallel class has size at most two, there are no other elements of $E(M)$ contained in $N_i$ alone, and so $b = a+1$. Now, assume that $i \neq 1$ and $i \neq r$. 
  
    Since $M$ is an excluded minor for the class of bicircular matroids, \autoref{lem:parallel_vert_3} implies that $M$ is vertically $3$-connected. Furthermore, since $a$ is not an element of $N_{i-1}$, we have $u_{i-1} < a$. Similarly, $b$ is not an element of $N_{i+1}$, so $l_{i+1} > b = a + 1$. Hence, $l_{i+1} > u_{i-1}$, violating the upper bound property and contradicting \autoref{lem:vert_3_connected_interval}. Therefore, we have $i = 1$ or $i = r$. The element $l_1$ is only contained in $N_1$, so if $i=1$ then we must have $a = l_1, b=l_1 + 1$; otherwise, $l_1$ is in parallel with $a$ and $b$, a contradiction. Similarly, if $i=r$ we must have $a = u_r - 1, b = u_r$. 
    Finally, the interval $[l_1, l_2 - 1]$ always consists of elements only in $N_1$, so we have a size two parallel class $\{l_1, l_1 + 1\}$ when $l_2 = l_1 + 2$ and otherwise $l_2 = l_1 + 1$. Similarly, $u_r = u_{r-1} + 2$ when $\{u_{r-1}, u_r\}$ is a parallel class, and $u_r = u_{r-1} + 1$ otherwise.
  \end{proof}
  
\begin{figure}[H]
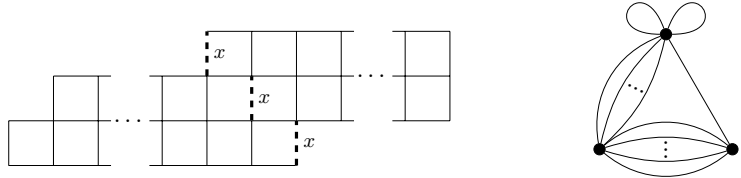

      \centering
      \begin{subfigure}[c]{0.7\textwidth}
          \centering
          \resizebox{0.85\width}{!}{
            \input{Tikz/MaximalGuys/Rank3/latticepath1}
          }
      \end{subfigure}
      \hfill
      \begin{subfigure}[c]{0.28\textwidth}
          \vspace*{-0.6cm}
          \resizebox{0.85\width}{!}{
            \input{Tikz/MaximalGuys/Rank3/bicircular1}
          }
      \end{subfigure}
      \caption{Lattice path and bicircular representations of the matroid constructed in Case 1.}
      \label{fig:no_more_rank_3_case_1}
  \end{figure}

\begin{figure}[H]
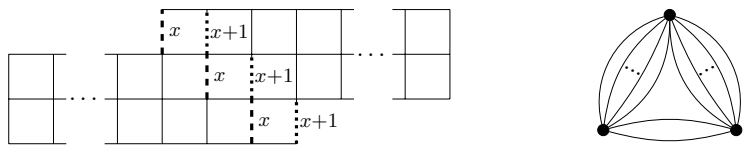

      \centering
      \begin{subfigure}[c]{0.7\textwidth}
          \centering
          \resizebox{0.85\width}{!}{
            \input{Tikz/MaximalGuys/Rank3/latticepath2}
          }
      \end{subfigure}
      \hfill
      \begin{subfigure}[c]{0.28\textwidth}
          \resizebox{0.85\width}{!}{
            \input{Tikz/MaximalGuys/Rank3/bicircular2}
          }
      \end{subfigure}
      \caption{Lattice path and bicircular representations of the matroid constructed in Case 2.}
      \label{fig:no_more_rank_3_case_2}
  \end{figure}

  \begin{lemma}
    Let $M$ be a rank-$3$ lattice path matroid. Then $M$ is an excluded minor for the class of bicircular and lattice path matroids if and only if $M$ is isomorphic to one of the matroids $U_{3,7}$, $T_3(U_{1,2} \oplus U_{3,5})$ and $T_3(U_{1,2} \oplus U_{1,2} \oplus U_{3,3})$.
    \label{lem:no_more_rank_3}
  \end{lemma}
  \begin{proof}
      Let $M$ be a rank-$3$ lattice path excluded minor for the class of bicircular and lattice path matroids and suppose that $M$ has the standard presentation $\mathcal{N} = (N_1,N_2,N_3)$ with $N_i = [l_i,u_i]$ for each $i \in \{1,2,3\}$. It follows from \autoref{lem:lattice_path_excluded_minors} that we here assume $M$ is not isomorphic to $U_{3,7}$, $T_3(U_{1,2} \oplus U_{3,5})$ or $T_3(U_{1,2} \oplus U_{1,2} \oplus U_{3,3})$. Since $M$ does not have a $U_{3,7}$-minor, there are at most two elements in $E(M)$ contained in all three of $N_1,N_2$ and $N_3$. Furthermore, $M$ is not bicircular, so \autoref{lem:2x3_subgrid} implies that there is at least one element in all three intervals. This gives us two cases:
  
      Case 1: Assume only one element, $x$, is an element of all three intervals of $\mathcal{N}$. Then $N_1$, $N_2$ and $N_3$ contain the intervals $[x-2, x]$, $[x-1, x+1]$ and $[x,x+2]$ respectively. Since $M$ is an excluded minor for the class of bicircular matroids, it follows from Lemmas \ref{lem:parallel_vert_3} and \ref{lem:uniform_minor} that $M$ has a $U_{3,5}$-minor. Therefore, if $M$ has more than one non-trivial parallel class, $M$ has a $T_3(U_{1,2} \oplus U_{1,2} \oplus U_{3,3})$-minor. Hence, $M$ has at most one non-trivial parallel class. Applying \autoref{lem:gross_parallel_class}, we can assume without loss of generality that this class is $\{l_1, l_1 + 1\}$, and therefore $l_2 = l_1 + 2$ and $u_3 = u_2 + 1$. It is easy to check that $M$ is a minor of the matroid $M[\mathcal{N}']$, given in \autoref{fig:no_more_rank_3_case_1} (left), where $\mathcal{N}' = (N_1',N_2',N_3')$ for $N_1' = [x-m-1,x]$, $N_2' = [x-m+1,x+n-1]$ and $N_3' = [x,x+n]$ for some non-negative integers $n$ and $m$. However, this matroid is bicircular, with the bicircular representation in \autoref{fig:no_more_rank_3_case_1} (right), so it is not an excluded minor for the class of bicircular and lattice path matroids.
  
      Case 2: Assume exactly two elements, $x$ and $x+1$, are elements of all three intervals of $\mathcal{N}$. Then $N_1$, $N_2$ and $N_3$ contain the intervals $[x-2, x+1]$, $[x-1, x+2]$ and $[x,x+3]$ respectively. If $M$ has a non-trivial parallel class then $T_3(U_{1,2} \oplus U_{3,5})$ is a minor of $M$. Hence, $M$ has no non-trivial parallel classes, and so $u_3 = u_2 + 1$ and $l_2 = l_1 + 1$ by \autoref{lem:gross_parallel_class}. It is easy to check that $M$ is a minor of the lattice path matroid $M[\mathcal{N}']$, given in \autoref{fig:no_more_rank_3_case_2} (left), where $\mathcal{N}' = (N_1', N_2', N_3')$ for $N_1' = [x-m,x+1]$, $N_2' = [x-m+1, x+n-1]$, and $N_3' = [x, x+n]$ with some non-negative integers $n$ and $m$. However, this matroid is also bicircular, with the bicircular representation in \autoref{fig:no_more_rank_3_case_2} (right), so it is not an excluded minor for the class of bicircular and lattice path matroids either.
      
      Hence, $M$ does not exist, and so our list of rank-$3$ lattice path excluded minors for the class of bicircular and lattice path matroids is complete.
  \end{proof}

\begin{figure}[H]
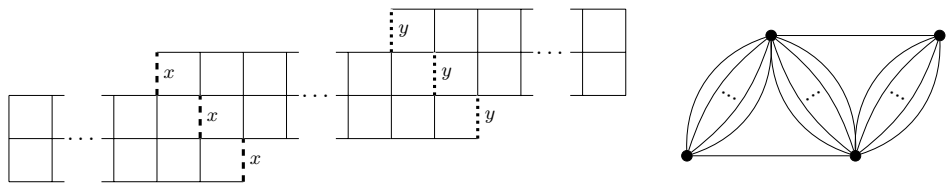

    \centering
    \hspace*{-1.2cm}
    \begin{subfigure}[c]{0.76\textwidth}
        \centering
        \resizebox{0.82\width}{!}{
          \input{Tikz/MaximalGuys/Rank4/latticepath2}
        }
    \end{subfigure}
    \begin{subfigure}[c]{0.22\textwidth}
        \centering
        \resizebox{0.8\width}{!}{
          \input{Tikz/MaximalGuys/Rank4/bicircular2}
        }
    \end{subfigure} 
    \caption{Lattice path and bicircular representations of $M[\mathcal{N''}]$.}
    \label{fig:no_more_rank_4}
  \end{figure}

  \begin{lemma}
    Let $M$ be a rank-$4$ lattice path matroid. Then $M$ is an excluded minor for the class of bicircular and lattice path matroids if and only if $M$ is isomorphic to one of the matroids $U_{4,7}$, 
    $T_4(U_{3,4} \oplus U_{3,3})$ and $T_4(U_{1,2} \oplus U_{4,5})$.
    \label{lem:no_more_rank_4}
  \end{lemma} 
  \begin{proof}
  Let $M$ be a rank-$4$ lattice path excluded minor for the class of bicircular and lattice path matroids, and suppose $M$ has the standard presentation $\mathcal{N} = (N_1,N_2,N_3,N_4)$ with $N_i = [l_i,u_i]$ for each $i \in \{1,2,3,4\}$. It follows from \autoref{lem:lattice_path_excluded_minors} that we here assume $M$ is not isomorphic to $U_{4,7}$, $T_4(U_{3,4} \oplus U_{3,3})$ or $T_4(U_{1,2} \oplus U_{4,5})$. 

  Since $M$ is not bicircular, \autoref{lem:2x3_subgrid} implies that $M$ has at least one element contained in three consecutive intervals of $\mathcal{N}$. Without loss of generality, let these intervals be $N_1,N_2$ and $N_3$. Assume $x$ is the largest such element. Furthermore, if the greatest element of $N_2$, $u_2$, is not an element of $N_4$, then $l_4 > u_2$, which violates the upper bound property, contradicting \autoref{lem:vert_3_connected_interval}. Thus, $u_2$ is contained in $N_2$, $N_3$ and $N_4$, and so $M$ has some least element $y$ contained in the intervals $N_2, N_3$ and $N_4$. Since $M$ has no $U_{4,7}$-minor, there are no elements contained in all four intervals of $\mathcal{N}$, so we have $y > x$.

  It follows that the intervals $N_1$, $N_2$, $N_3$ and $N_4$ must contain the sub-intervals $N_1' = [x-2,x]$, $N_2' = [x-1,y]$, $N_3' = [x,y+1]$ and $N_4' = [y,y+2]$ respectively.
  The lattice path matroid $M' = M[(N_1',N_2',N_3',N_4')]$ is vertically $3$-connected, so \autoref{lem:uniform_minor} implies that $M'$ has a $U_{4,6}$-minor. Hence, if $M$ has a non-trivial parallel class, $M$ has a $T_4(U_{1,2} \oplus U_{4,5})$-minor, a contradiction. Thus, $M$ has no non-trivial parallel classes, and so all elements of $E(M)$ are contained in at least two intervals of $\mathcal{N}$. Similarly, if $M$ has an element other than $x$ or $y$ contained in exactly three intervals of $\mathcal{N}$ then $M$ has a $T_4(U_{3,4} \oplus U_{3,3})$-minor, another contradiction. Hence, the only elements contained in exactly three intervals of $\mathcal{N}$ are (at most) $x$ and $y$. In summary, all elements of $M$ other than $x$ and $y$ must be contained in at most two intervals of $\mathcal{N}$, and since $M$ has no non-trivial parallel classes, only $l_1$ and $u_4$ are contained in exactly one interval of $\cN$.

  It is therefore easy to see that $M$ is a minor of the matroid $M[\mathcal{N''}]$ given in \autoref{fig:no_more_rank_4} (left), where $\mathcal{N''}$ consists of the intervals $N_1'' = [x-a,x]$, $N_2'' = [x-a+1,y]$, $N_3'' = [x,y+b-1]$ and $N_4'' = [y,y+b]$ for $a \geq 2$ and $b \geq 2$ and some $y > x$. A  bicircular representation of this matroid is given in \autoref{fig:no_more_rank_4} (right), and so $M$ is bicircular. Hence, $M$ is not an excluded minor for the class of bicircular and lattice path matroids.
\end{proof}

  By \autoref{lem:2x3_subgrid}, all lattice path matroids with rank at most two are bicircular. Thus, by combining Lemmas \ref{lem:only_u_57}, \ref{lem:no_more_rank_3} and \ref{lem:no_more_rank_4}, we establish that the matroids in part (i) of \autoref{thm:main} form the complete list of lattice path excluded minors for the class of bicircular and lattice path matroids.

\section{Excluded Minors that are Not Lattice Path} \label{sec:non_lattice_path_EMs}

In this section, we determine the list of non-lattice path excluded minors for the class of bicircular and lattice path matroids, thereby completing the proof of \autoref{thm:main}. 

We begin by noting that, by \autoref{thm:lattice_path_excluded_minors}, each of the matroids $A_3$, $B_{3,3}$, $C_{4,2}$, $C_{5,2}$, $D_4$, $\mathcal{W}^3$, $R_3$, $R_4$, $B_{2,2}$, $B_{3,2}$, $E_4$ and $\mathcal{W}_3$ is an excluded minor for the class of lattice path matroids. Moreover, it is easily checked that each of the matroids $A_3$, $B_{3,3}$, $C_{4,2}$, $C_{5,2}$, $D_4$, $\mathcal{W}^3$, $R_3$ and $R_4$ are bicircular by checking the bicircular representations given in \autoref{fig:excluded_minors_main_2}. Finally, by \autoref{lem:bicircular_excluded_minors}, each of the matroids $B_{2,2}$, $B_{3,2}$, $E_4$ and $\mathcal{W}_3$ are excluded minors for the class of bicircular matroids. This gives us the following result:

\begin{lemma}
  The bicircular matroids $A_3$, $B_{3,3}$, $C_{4,2}$, $C_{5,2}$, $D_4$, $\mathcal{W}^3$, $R_3$ and $R_4$ are excluded minors for the class of bicircular and lattice path matroids.
  The non-bicircular and non-lattice path matroids $B_{2,2}$, $B_{3,2}$, $E_4$ and $\mathcal{W}_3$ are excluded minors for the class of bicircular and lattice path matroids.
  \label{lem:non_lattice_excluded_b_n_l}
\end{lemma}

To prove that \autoref{lem:non_lattice_excluded_b_n_l} gives the complete lists of bicircular and neither bicircular nor lattice path excluded minors for the class of bicircular and lattice path matroids, it suffices to show that each of the remaining excluded minors for the class of lattice path matroids has a proper minor that is not bicircular.

We first provide a bound on the number of elements in a bicircular matroid.
\begin{lemma}
  Let $M$ be a rank-$r$ bicircular matroid. If every circuit of $M$ has more than $k$ elements for some $k > 2$, then $E(M) \leq \frac{r(r-1)}{k-2}$.
  \label{lem:bicircular_bound_1}
\end{lemma}
\proof Let $M$ be a rank-$r$ bicircular matroid. Then \autoref{lem:bicircular_simplex} implies that 
every element of $M$ can be placed on the edges or vertices of an $r$-vertex simplex $\Delta$. 
If every circuit of $M$ has more than $k$ elements, then any $(k-1)$-face of $\Delta$ contains a maximum of $k-1$ elements. 
There are $r \choose k-1$ total $(k-1)$-faces of $\Delta$, and every element of $M$ lies on a minimum of $r-2 \choose k-3$ $(k-1)$-faces of $\Delta$. 
Hence, the total number of elements in $M$ is bounded above by 
\begin{align*}
  |E(M)| &\leq \frac{(k-1){r \choose k-1}}{{r-2 \choose k-3}}\\
        &= \frac{r(r-1)}{k-2}
\end{align*}\qed 

We now check the members of each of the infinite families of matroids given as excluded minors for the class of lattice path matroids in \autoref{thm:lattice_path_excluded_minors}.
\begin{lemma}
  For all $n \geq 3$, the matroid $A_n$ is an excluded minor for the class of bicircular and lattice path matroids if and only if $n = 3$.
  \label{lem:an_excluded_minors}
\end{lemma}
  
\begin{proof}
  If $n=3$, then, by \autoref{lem:non_lattice_excluded_b_n_l}, $A_3$ is an excluded minor for the class of bicircular and lattice path matroids. Let $n \geq 4$, and let $f$ be the unique element in the intersection of the two distinct $n$-element circuits of $A_n$. Consider $A_n \backslash f$. By construction, it is easily seen that $A_n \backslash f$ is isomorphic to $U_{n,2n-1}$. Therefore, if $n \geq 4$, then $A_n$ has a proper minor isomorphic to $U_{4,7}$. Since $U_{4,7}$ is an excluded minor for the class of bicircular and lattice path matroids, it follows that $A_n$ is not an excluded minor for this class. This completes the proof of the lemma.
\end{proof}

\begin{lemma}
  For all $n \geq k \geq 2$, the matroid $B_{n,k}$ is an excluded minor for the class of bicircular and lattice path matroids if and only if $(n,k) \in \{(2,2),(3,2),(3,3)\}$.
\end{lemma}
\begin{proof}
  
If $(n,k) \in \{(2,2),(3,2),(3,3)\}$, then, by \autoref{lem:non_lattice_excluded_b_n_l}, $B_{2,2}$, $B_{3,2}$ and $B_{3,3}$ are excluded minors for the class of bicircular and lattice path matroids. Hence, suppose $B_{n,k}$ is an excluded minor for the class of bicircular and lattice path matroids with $n > 3$.

Let $f$ be an element of the $k$-element circuit of $B_{n,k}$, and let $g$ and $h$ be elements from each of the $n$-element circuits of $B_{n,k}$. Consider $B_{n,k}\backslash \{f,g,h\}$. This is a rank-$n$ matroid on $2n + k - 3$ elements, in which every circuit has more than $n$ elements. Therefore, since $n > 3$, \autoref{lem:bicircular_bound_1} implies that if
$B_{n,k}\backslash \{f,g,h\}$ is bicircular, then
$$2n+k - 3 \leq \frac{n(n-1)}{n-2},$$ which in turn implies
$$n \leq \frac{6-k + \sqrt{k^2 -4k + 12}}{2}.$$

The sequence $\{x_k\}_{k=2}^\infty = \left\{\frac{6-k + \sqrt{k^2 -4k + 12}}{2}\right\}_{k=2}^\infty$ is monotonically decreasing, and we have, by definition of $B_{n,k}$, that $n \geq k \geq 2$, so this leaves only a few cases to check.
If $k=2$ then we have $x_2 = \frac{4 + \sqrt{8}}{2}$ and so $n \leq 3$. If $k=3$ then $x_3 = 3$ and so $n=3$. 
No larger integer values of $k$ satisfy $n \geq k$, so $B_{n,k}$ has non-bicircular proper minors unless $B_{n,k}$ is $B_{3,2}$, $B_{2,2}$ or $B_{3,3}$. This completes the proof of the lemma.
\end{proof}

\begin{lemma}
  For all $n \geq k \geq 2$, the matroid $C_{n+k,k}$ is an excluded minor for the class of bicircular and lattice path matroids if and only if $(n+k,k) \in \{(4,2), (5,2)\}$.
\end{lemma}
  
\begin{proof}
  Suppose that $C_{n+k,k}$ is an excluded minor for the class of bicircular matroids, and consider its dual, $B_{n,k}$. Let $f$ be an element of the $k$-element circuit $C$ of $B_{n,k}$, and consider $B_{n,k}/f$. By construction, it is easily seen that $B_{n,k}/f$ is isomorphic to $T_n(U_{n-1,2n} \oplus U_{k-2,k-1})$. Thus, if $k=2$, then $B_{n,k}/f\backslash (C-f)$ is isomorphic to $U_{n-1,2n}$, while, if $k \geq 3$, then $B_{n,k}/f \backslash (C-\{f,g\})$, where $g \in (C - f)$, is isomorphic to $U_{n-1,2n+1}$. In particular, if $k=2$, then $C_{n+k,k}$ has a proper $U_{n+1,2n}$-minor and, if $k \geq 3$, then $C_{n+k,k}$ has a proper $U_{n+2,2n+1}$-minor. In turn, this implies that if $k=2$ and $n \geq 4$, or if $n \geq k \geq 3$, then $C_{n+k,k}$ has $U_{5,7}$ as a proper minor. Therefore, as $U_{5,7}$ is an excluded minor for the class of bicircular and lattice path matroids, it follows that if $k=2$ and $n \geq 4$, or if $n \geq k \geq 3$, then $C_{n+k,k}$ is not an excluded minor for this class.

  Conversely, if $(n+k,k) \in \{(4,2), (5,2)\}$, then \autoref{lem:non_lattice_excluded_b_n_l} implies that $C_{n+k,k}$ is an excluded minor for the class of bicircular and lattice path matroids. This completes the proof of the lemma.
\end{proof}
  
\begin{lemma}
  For all $n \geq 4$, the matroids $D_n$ and $E_n$ are excluded minors for the class of bicircular and lattice path matroids if and only if $n = 4$.
  \label{lem:dn_en_excluded_minors}
\end{lemma}
  
\begin{proof} 
  If $n=4$, then $D_4$ and $E_4$ are excluded minors for the class of bicircular and lattice path matroids by \autoref{lem:non_lattice_excluded_b_n_l}. Let $n \geq 5$, let $\{x,y\}$ be the unique series pair of $D_n$, and let $C_1$ and $C_2$ be the disjoint $(n-1)$-element circuits of $D_n \backslash \{x,y\}$. If $X$ is a $(n-1)$-element subset of $C_1 \cup C_2$ that is neither $C_1$ nor $C_2$, then $X \cup \{x,y\}$ is a $(n+1)$-element circuit of $D_n$. Hence, every $(n-1)$-element subset of $D_n/\{x,y\}$ is a circuit and so $D_n/\{x,y\}$ is isomorphic to $U_{n-2,2n-2}$. Thus, $D_n$ has a $U_{n-2,2n-2}$-minor. It now follows that, since $n \geq 5$, then $D_n$ has a proper $U_{3,8}$-minor, and hence a proper $U_{3,7}$-minor. Since $U_{3,7}$ is an excluded minor for the class of bicircular matroids, this implies that $D_n$ is not an excluded minor for the class of bicircular and lattice path matroids.
  
  Finally, since $D_n$ has a proper $U_{3,8}$-minor, it also follows that $D_n$ has a proper $U_{2,7}$-minor. Therefore, since $E_n$ is the dual of $D_n$, the matroid $E_n$ has a proper $U_{5,7}$-minor. Since $U_{5,7}$ is also an excluded minor for the class of bicircular matroids, it follows that $E_n$ is not an excluded minor for the class of bicircular and lattice path matroids for all $n \geq 5$.
\end{proof}

Hence, by Lemmas \ref{lem:an_excluded_minors}--\ref{lem:dn_en_excluded_minors}, the excluded minors for the class of lattice path matroids that are not listed in \autoref{lem:non_lattice_excluded_b_n_l} all have minors that are not bicircular. Therefore, we have established that the matroids in part (ii) of \autoref{thm:main} form the complete list of bicircular excluded minors for the class of bicircular and lattice path matroids, and that the matroids listed in part (iii) of \autoref{thm:main} form the complete list of excluded minors for the class of bicircular and lattice path matroids that are neither bicircular nor lattice path. This completes the proof of \autoref{thm:main}.

\bibliographystyle{elsarticle-num}
\bibliography{references}

\end{document}